\documentclass{ws-idaqp}

\usepackage{ifthen}
\newboolean{article}
    \setboolean{article}{true}
\newboolean{report}
\newboolean{book}
\newboolean{letter}
\newboolean{german}
\newboolean{italian}
\newboolean{nobaselinestretch}
   \setboolean{nobaselinestretch}{true}
\newboolean{nosectionappendix}
\newboolean{oldtoc}
\newboolean{nosectionequn}
\newboolean{notheorem}

\ifthenelse{\boolean{german}}{
    \usepackage{german}}{}

\usepackage[latin1]{inputenc}

\usepackage{amsmath}
\usepackage{amssymb}
\usepackage[mathscr]{eucal}

\ifthenelse{\boolean{notheorem}}{}{
    \usepackage{theorem}}

\ifthenelse{\boolean{nobaselinestretch}}{}{
    \renewcommand{\baselinestretch}{1.25}}

\newenvironment{env}[2]{\begin{#1}#2\end{#1}}{}
    \newcommand{\beq}[1]{\begin{env}{equation}{#1}}
    \newcommand{\beqn}[1]{\begin{env}{equation*}{#1}}
    \newcommand{\bal}[1]{\begin{env}{align}{#1}}
    \newcommand{\baln}[1]{\begin{env}{align*}{#1}}
    \newcommand{\bga}[1]{\begin{env}{gather}{#1}}
    \newcommand{\bgan}[1]{\begin{env}{gather*}{#1}}
    \newcommand{\bflal}[1]{\begin{env}{flalign}{#1}}
    \newcommand{\bflaln}[1]{\begin{env}{flalign*}{#1}}
    \newcommand{\bmu}[1]{\begin{env}{multline}{#1}}
    \newcommand{\bmun}[1]{\begin{env}{multline*}{#1}}
    \newcommand{\bsp}[1]{\begin{env}{split}{#1}}

    \newcommand{\eeq}{\end{env}}
    \newcommand{\eeqn}{\end{env}}
    \newcommand{\eal}{\end{env}}
    \newcommand{\ealn}{\end{env}}
    \newcommand{\ega}{\end{env}}
    \newcommand{\egan}{\end{env}}
    \newcommand{\eflal}{\end{env}}
    \newcommand{\eflaln}{\end{env}}
    \newcommand{\emu}{\end{env}}
    \newcommand{\emun}{\end{env}}
    \newcommand{\esp}{\end{env}}

\newcommand{\lf}{\vspace{2ex}}
\newcommand{\bulletline}[1][]{\lf\noindent~\hfill$\bullet\bullet\bullet$\hfill~

\lf\noindent\bf{#1}}

\renewcommand{\bf}[1]{\textbf{#1}}
\renewcommand{\it}[1]{\textit{#1}}
\renewcommand{\sc}[1]{\textsc{#1}}
\renewcommand{\sf}[1]{\textsf{#1}}
\renewcommand{\sl}[1]{\textsl{#1}}
\renewcommand{\tt}[1]{\texttt{#1}}

\newcommand{\hl}[1]{\bf{\it{#1}}}
\newcommand{\mrm}[1]{\mathrm{#1}}
\newcommand{\mbf}[1]{\mathbf{#1}}
\newcommand{\msf}[1]{\text{\small$\sf{#1}$}}

\newcommand{\cmc}[1]{\mathcal{#1}}
\newcommand{\eus}[1]{\mathscr{#1}}

\newcommand{\bb}[1]{\mathbb{#1}}
\newcommand{\msmall}[1]{{\setlength{\arraycolsep}{.6ex}\text{\small$#1$}}}

\newcommand{\mtiny}[1]{{\setlength{\arraycolsep}{.3ex}\text{\tiny$#1$}}}
\newcommand{\nbd}[1]{$#1$\nobreakdash--}
\newcommand{\ol}[1]{\overline{#1}}

\newcommand{\vp}{\varphi}

\newcommand{\abs}[1]{\left\lvert#1\right\rvert}
\newcommand{\norm}[1]{\left\lVert#1\right\rVert}

\newcommand{\bfam}[1]{\bigl(#1\bigr)}
\newcommand{\Bfam}[1]{\Bigl(#1\Bigr)}
\newcommand{\AB}[1]{\langle#1\rangle}
\newcommand{\bAB}[1]{\bigl\langle#1\bigr\rangle}
\newcommand{\BAB}[1]{\Bigl\langle#1\Bigr\rangle}
\newcommand{\CB}[1]{\{#1\}}
\newcommand{\bCB}[1]{\bigl\{#1\bigr\}}
\newcommand{\BCB}[1]{\Bigl\{#1\Bigr\}}
\newcommand{\SB}[1]{[#1]}

\newcommand{\LO}[1]{(#1]}

\newcommand{\Matrix}[1]{\begin{pmatrix}#1\end{pmatrix}}
\newcommand{\SMatrix}[1]{\msmall{\Matrix{#1}}}

\newcommand{\tMatrix}[1]{\mtiny{\Matrix{#1}}}
\newcommand{\rtMatrix}[1]{\raisebox{.3ex}{\tMatrix{#1}}}

\newcommand{\set}[2][]{
    \ifthenelse{\equal{#1}{}}{
        \CB{#2}}{
        \CB{#1~|~#2}}}
\newcommand{\bset}[2][]{
    \ifthenelse{\equal{#1}{}}{
        \bCB{#2}}{
        \bCB{#1~|~#2}}}
\newcommand{\Bset}[2][]{
    \ifthenelse{\equal{#1}{}}{
        \BCB{#2}}{
        \BCB{#1~\big|~#2}}}
\newcommand{\zero}{\CB{0}}

\DeclareMathOperator{\id}{\normalfont\msf{id}}

\renewcommand{\ker}{\operatorname{\msf{ker}}}
\renewcommand{\dim}{\operatorname{\msf{dim}}}

\newcommand{\C}{\bb{C}}

\newcommand{\E}{\bb{E}}

\newcommand{\R}{\bb{R}}

\newcommand{\cB}{\cmc{B}}

\newcommand{\sB}{\eus{B}}

\newcommand{\sS}{\eus{S}}

\newcommand{\U}{\mbf{1}}
\newcommand{\F}{{\mrm{F}}}

\ifthenelse{\boolean{nosectionequn}}{}{
    \numberwithin{equation}{section}
    \renewcommand{\theequation}{\thesection.\arabic{equation}}}

\ifthenelse{\boolean{article}\or\boolean{letter}\or\boolean{nosectionequn}}{
    \setboolean{nosectionappendix}{true}}{}
\ifthenelse{\boolean{nosectionappendix}}{}{
    \renewcommand{\appendix}{
        \chapter*{\appendixname}
        \addcontentsline{toc}{chapter}{\appendixname}
        \renewcommand{\thesection}{\Alph{section}}
        \setcounter{section}{0}}}
   
\ifthenelse{\boolean{report}\or\boolean{book}}{
    }{}

\ifthenelse{\boolean{notheorem}}{}{
        \newcommand{\notename}{Note.}
        \newcommand{\mnname}{Mathematical note.}
        \newcommand{\enname}{End of the note.}
        \newcommand{\definame}{Definition.}
        \newcommand{\propname}{Proposition.}
        \newcommand{\lemname}{Lemma.}
        \newcommand{\exname}{Example.}
        \newcommand{\exername}{Exercise.}
        \newcommand{\remname}{Remark.}
        \newcommand{\obname}{Observation.}
        \newcommand{\thmname}{Theorem.}
        \newcommand{\corname}{Corollary.}
        \newcommand{\proofname}{Proof.}
    \ifthenelse{\boolean{german}}{
        \renewcommand{\mnname}{Mathematische Notiz.}
        \renewcommand{\enname}{Ende der Notiz.}
        \renewcommand{\exname}{Beispiel.}
        \renewcommand{\exername}{Übung.}
        \renewcommand{\remname}{Bemerkung.}
        \renewcommand{\obname}{Beobachtung.}
        \renewcommand{\thmname}{Satz.}
        \renewcommand{\corname}{Korollar.}
        \renewcommand{\proofname}{Beweis.}}{}
    \ifthenelse{\boolean{italian}}{
        \renewcommand{\mnname}{Nota matematica.}
        \renewcommand{\enname}{Fina della nota.}
        \renewcommand{\definame}{Definizione.}
        \renewcommand{\propname}{Proposizione.}
        \renewcommand{\exname}{Esempio.}
        \renewcommand{\exername}{Esercizio.}
        \renewcommand{\remname}{Nota.}
        \renewcommand{\obname}{Osservazione.}
        \renewcommand{\thmname}{Teorema.}
        \renewcommand{\corname}{Corollario.}
        \renewcommand{\proofname}{Dimostrazione.}

       \renewcommand{\bibname}{Bibliografia}
       
       \renewcommand{\appendixname}{Appendice}
       \renewcommand{\contentsname}{Indice}
       \renewcommand{\listfigurename}{Elenco delle figure}
       \renewcommand{\listtablename}{Elenco delle tabelle}
       
       \renewcommand{\figurename}{Figura}
       \renewcommand{\tablename}{Tabella}

       }{}
    \theoremheaderfont{\normalfont\bfseries}
    \theoremstyle{change}
        \theorembodyfont{\rmfamily}
            \newtheorem{emp}{}[section]
                \newcommand{\bemp}[1][]{
                    \begin{emp}\hskip-\labelsep\bf{#1}\hskip\labelsep}
                \newcommand{\eemp}{\end{emp}}
\newtheorem{itemp}[emp]{}
                \newcommand{\bitemp}[1][]{
                    \begin{itemp}\hskip-\labelsep\bf{#1}\hskip\labelsep\normalfont\itshape}
                \newcommand{\eitemp}{\end{itemp}}
            \newtheorem{note}[emp]{\notename}
                \newcommand{\bnote}{\begin{note}}
                \newcommand{\enote}{\end{note}}
            \newtheorem{mn}[emp]{\mnname}
                \newcommand{\bnm}{\begin{mn}~\begin{quotation}\renewcommand{\baselinestretch}{1}\small\noindent\ignorespaces}
                \newcommand{\enm}{\end{quotation}\hfill\bf{\enname}\end{mn}}
            \newtheorem{ex}[emp]{\exname}
                \newcommand{\bex}{\begin{ex}}
                \newcommand{\eex}{\end{ex}}
            \newtheorem{exer}[emp]{\exername}
                \newcommand{\bexer}{\begin{exer}}
                \newcommand{\eexer}{\end{exer}}
            \newtheorem{defi}[emp]{\definame}
                \newcommand{\bdefi}{\begin{defi}}
                \newcommand{\edefi}{\end{defi}}
            \newtheorem{rem}[emp]{\remname}
                \newcommand{\brem}{\begin{rem}}
                \newcommand{\erem}{\end{rem}}
            \newtheorem{ob}[emp]{\obname}
                \newcommand{\bob}{\begin{ob}}
                \newcommand{\eob}{\end{ob}}
        \theorembodyfont{\normalfont\itshape}
            \newtheorem{thm}[emp]{\thmname}
                \newcommand{\bthm}{\begin{thm}}
                \newcommand{\ethm}{\end{thm}}
            \newtheorem{prop}[emp]{\propname}
                \newcommand{\bprop}{\begin{prop}}
                \newcommand{\eprop}{\end{prop}}
            \newtheorem{cor}[emp]{\corname}
                \newcommand{\bcor}{\begin{cor}}
                \newcommand{\ecor}{\end{cor}}
            \newtheorem{lem}[emp]{\lemname}
                \newcommand{\blem}{\begin{lem}}
                \newcommand{\elem}{\end{lem}}
\newenvironment{empn}[1]{\lf\noindent\bf{#1}\ignorespaces\hskip\labelsep}{\lf}
		\newcommand{\bempn}[1]{\begin{empn}{#1}}
		\newcommand{\eempn}{\end{empn}}
		\newcommand{\bitempn}[1]{\begin{empn}{#1}\normalfont\itshape}
		\newcommand{\eitempn}{\end{empn}}
                \newcommand{\bnmn}{\begin{empn}{\mnname}~\begin{quotation}\renewcommand{\baselinestretch}{1}\small\noindent\ignorespaces}
                \newcommand{\enmn}{\end{quotation}\hfill\bf{\enname}\end{empn}}
		\newcommand{\bexn}{\begin{empn}{\exname}}
		\newcommand{\eexn}{\end{empn}}
		\newcommand{\bexern}{\begin{empn}{\exername}}
		\newcommand{\eexern}{\end{empn}}
		\newcommand{\bdefin}{\begin{empn}{\definame}}
		\newcommand{\edefin}{\end{empn}}
		\newcommand{\bremn}{\begin{empn}{\remname}}
		\newcommand{\eremn}{\end{empn}}
		\newcommand{\bobn}{\begin{empn}{\obname}}
		\newcommand{\eobn}{\end{empn}}

\newcommand{\qedsymbol}{~\hfill$\oblong$}
    \newcounter{proof}[emp]
    \newenvironment{Proof}[1]{
        \vspace{1ex}
        \renewcommand{\item}[1][\stepcounter{proof}(\roman{proof})]%
            {##1\hskip\labelsep}
        \noindent\bf{#1\hskip\labelsep}}{
        \nolinebreak\qedsymbol}
     \renewcommand{\proof}[1][\proofname]{
        \begin{Proof}{#1}\ignorespaces}
    \renewcommand{\qed}{\end{Proof}}
    \newcommand{\noqed}{
        \renewcommand{\qedsymbol}{}
        \end{Proof}}}
    \ifthenelse{\boolean{italian}}{
        \renewcommand{\proofname}{Dimostrazione.}}{}

\newcounter{OP}
		\newcommand{\bOP}{\stepcounter{OP}\begin{empn}{Open Problem \arabic{OP}:}}
		\newcommand{\eOP}{\end{empn}}

\usepackage{lmodern}

\usepackage{stmaryrd}

\usepackage[matrix,arrow,curve]{xy}

\usepackage{hyperref}

\begin{document}

\markboth{Michael Skeide}{Partial Isometries}

\bibliographystyle{amsalpha}

\title{Partial Isometries Between Hilbert Modules\\and Their Compositions\\~\\\normalfont To the Memory of Kalyanapuram Rangachari Parthasarathy}


\author{Michael Skeide}

\address{Dipartimento di Economia, Universit\`a\ degli Studi del Molise, Via de Sanctis\\
86100 Campobasso, Italy\\
skeide@unimol.it
}

\maketitle

\begin{abstract}
Motivated by questions raised in the preprint \cite{AcLu20p} by Accardi and Lu (private communication), we examine criteria for when the product of two partial isometries between Hilbert spaces is again a partial isometry and we use this to define a new composition operation that always yields again a partial isometry. Then, we aim at promoting these results to (not necessarily adjointable) partial isometries between Hilbert modules as proposed by Shalit and Skeide \cite{ShaSk23}.

The case of Hilbert spaces is elementary and rather simple, though not trivial, but -- we expect -- folkloric. The case of Hilbert modules suffers substantially from the fact that bounded right linear maps need not possess necessarily an adjoint. In fact, we show that the new composition law for partial isometries between Hilbert spaces can in no way be promoted directly to partial isometries between Hilbert modules, but that we have to pass to the more flexible class of \it{partially defined isometries}.
\end{abstract}





\section{Introduction} \label{introSEC}

If $A,B,C,\ldots$ are sets, then a function $f\colon B\rightarrow A$ and function $g\colon C\rightarrow B$ can be composed to a function $f\circ g\colon C\rightarrow A$ defined as $f\circ g(c):=f(g(c))$. In the sequel, we look at functions, written as pairs $(f,D_f)$, that are only \hl{partially defined} as functions $f\colon D_f\rightarrow A$ on their domain $D_f\subset B$; we denote this situation as $f\colon B\supset D_f\rightarrow A$. In order to compose such partially defined functions $(f,D_f)$ and $(g,D_g)$, we have to indicate also a suitable domain $D_{f\circ g}$ for their composition:
\beq{ \label{pcdef}
(f,D_f)\circ(g,D_g)
~:=~
(f\circ g,D_{f\circ g})
}\eeq
where $D_{f\circ g}:=\bCB{c \in D_g\colon g(c)\in D_f}$ is the canonical maximal domain such that $f\circ g(c):=f(g(c))$ makes sense for $c\in D_{f\circ g}$.

\brem
In the preprint \cite{AcLu20p} by Accardi and Lu (private communication), the authors start with the partially defined functions and, then, aim at pushing their behaviour forward to certain partial isometries between interacting Fock spaces (the way the latter have been defined in Gerhold and Skeide \cite{GeSk20b}). In particular, their interest is turning interacting Fock spaces into a category with certain partial isometries as morphisms (adding to those in \cite{GeSk20b}, where only several classes of isomorphisms have been considered). Since in a category, by definition, for all objects $A,B,C$ and for all morphisms $f\colon B\rightarrow A$ and $g\colon C\rightarrow B$ there always exists the composition, a morphism $f\circ g\colon C\rightarrow A$, this raises the question (not answered in \cite{AcLu20p}) how to define a composition of two partial isometries $v$ and $w$ also in the case when their operator product $vw$ is no longer a partial isometry.

We, in these notes, are not really interested in categorical questions. But since we feel it has some relevance for the setting of \cite{AcLu20p}, we state what we can in a series of remarks -- remarks that, however, may safely be ignored.
\erem

\lf
To lift partially defined functions to partial isometries, it is indispensable to restrict to the \hl{injective partially defined functions} (in the sense that the function $f\colon D_f\rightarrow A$ is injective in the usual sense).

\brem \label{injrem}
Since the composition of injective partially defined functions is an injective partially defined function, the class $\sS$ with the injective partially defined functions is a subcategory of $\sS$ with the partially defined functions.
\erem

\lf
Now, if we associate with any set $A\in\sS$ the Hilbert space $H_A$ with orthonormal basis $\bfam{e_a}_{a\in A}$, then each injective partially defined function $f\colon B\supset D_f\rightarrow A$ gives rise to a partial isometry $v_f\colon H_B\rightarrow H_A$ defined by
\beqn{
v_f(e_b)
~:=~
\begin{cases}
e_{f(b)}&b\in D_f,
\\
0&\text{otherwise.}
\end{cases}
}\eeqn
(If $f$ is not injective, then $v_f$ would not be a partial isometry; not even a contraction. In the worst case, when a certain value $a\in A$ is taken by $f(b)$ for an infinite number of $b\in B$, then $v_f$ is not even well-defined.)

The following is easy to see.

\bob
Let $(f,D_f)$ be an injective partially defined function from $B\supset D_f$ to $A$. Then $v_f$ is an isometry if and only if $D_f=B$ (that is, if $f$ is everywhere defined).  $v_f$ is a coisometry if and only if $f$ is surjective. Consequently, $v_f$ is unitary if and only if (recall that $f$ is assumed injective!) $f$ is bijective from $B$ to $A$.
\eob

The following is almost as easy to see.

\bob
With $f$, $g$, and $f\circ g$ as for \eqref{pcdef}, we have
\beqn{
v_fv_g
~=~
v_{f\circ g}.
}\eeqn
In particular, the composition of these partial isometries is again a partial isometry. (The latter also follows from Lemma \ref{pilem} below, by simply observing that any pair of projections of the form $p=\sum_{b\in S}e_be_b^*$ for subsets $S$ of $B$ commute.)
\eob

\brem
We see that by $A\mapsto H_A$ we embed the class of sets into the class of Hilbert spaces, by $(f,D_f)\mapsto v_f$ we embed the class of injective partially defined functions into the class of partial isometries, and the latter embedding respects compositions for the usual product of partial isometries. However, not all partial isometries between Hilbert spaces have a product that is again a partial isometry, so the class of Hilbert spaces with partial isometries as morphisms is not a category under the usual operator composition.
\erem

\lf
In Section \ref{picSEC}, imitating the idea of giving the composition of partially defined functions its natural maximal domain, we provide a composition of partial isometries that is again a partial isometry. (This composition does turn Hilbert spaces with partial isometries into a category. And for Hilbert spaces $H_A$ and partial isometries $v_f$ we recover the usual composition.)

In Section \ref{modpiSEC} we perform the same program for Hilbert \nbd{\cB}modules ($\cB$ a fixed \nbd{C^*}algebra) and partial isometries between them. While for Hilbert spaces every bounded operator has an adjoint, for Hilbert modules this is not so. On the one hand, isometries need not have an adjoint. On the other hand, it is not a task with a too obvious solution to come up with a definition of coisometry that does not assume as part of the definition that the coisometry has an adjoint; the same is true for partial isometries. Here, we follow the definitions (without hypothesizing existence of an adjoint) from Shalit and Skeide \cite{ShaSk23}. (While isometries need not be adjointable, in \cite[Appendix B]{ShaSk23} it is shown that coisometries turn out to be always adjointable and are equivalently characterized by their usual definition; partial isometries need not be adjointable.) Most results of the Hilbert space case in Section \ref{picSEC} have to be reformulated before we can show them in Section \ref{modpiSEC} for Hilbert modules. Some proofs are simply different; some proofs work by reduction to the corresponding Hilbert space results via \bf{the tool}, a standard tool in Hilbert module theory (well known, but less widespread) we explain in the proof of  Theorem \ref{cMthm} and Footnote \ref{tpFN}; see also Remark \ref{SSlemrem} . The main new feature for modules is that to indicate a composition law for all (everywhere-defined!) partial isometries, we are forced to leave the class of partial isometries and have to replace it by the class of \it{partially defined isometries}, which are the perfect analogue of the \it{injective partially defined functions} among sets.
\vspace{-2ex}

\bulletline[A personal remark.~]
This little note is to commemorate our friend, colleague, teacher, advisor, mentor, ..., Partha, whom we are missing so much. While I did not have the fortune to have collaborated with Partha, I still had nice and instructive chats with him (about maths and not) or got ideas from listening to his talks  and lectures or reading his papers; a certain amount of my own work can be traced back to things he wrote and said. I think he enjoyed asking questions that can be formulated very simply, and finding answers to them. I hope this little note can meet at least a bit this style of his.

\newpage

\section{About composing partial isometries between Hilbert spaces} \label{picSEC}

Recall that an operator $p$ on a Hilbert space is a \hl{projection} if $p^*p=p$. (Equivalently, $p=p^*$ and $p^2=p$.) Recall that an operator $v$ between two Hilbert spaces is a \hl{partial isometry} if $vv^*v=v$. (Equivalently, $v^*v$ is a projection.%
\footnote{
Indeed, if $v^*v$ is a projection, then
\bmun{
~~~~~~
(vv^*v-v)^*(vv^*v-v)
~=~
(v^*vv^*-v^*)(vv^*v-v)
~=~
v^*v-v^*v-v^*v+v^*v
~=~
0,~~~~~~
}\emun
so $vv^*v-v=0$. Of course, if $v$ is a partial isometry, $(v^*v)^*v^*v=v^*vv^*v=v^*v$.
}%
) This happens if and only if $v^*$ is a partial isometry. (Simply, take adjoints.) So, also $vv^*$ is a projection.%

\blem \label{pilem}
For two partial isometries $v$ and $w$ the following are equivalent:
\begin{enumerate}
\item \label{pi1}
$vw$ is a partial isometry.

\item \label{pi2}
$v^*vww^*$ is an idempotent.

\item \label{pi3}
$v^*vww^*$ is a projection.

\item \label{pi4}
The projections $v^*v$ and $ww^*$ commute.
\end{enumerate}
\elem

\proof
Clearly, \ref{pi4} $\Rightarrow$ \ref{pi3} $\Rightarrow$ \ref{pi2}.

By definition, $vw$ is a partial isometry if and only if
\beq{ \label{vwpi}
vww^*v^*vw
~=~
vw.
}\eeq
Multiplying \eqref{vwpi} from the left with $v^*$ and from the right with $w^*$ gives
\beq{ \label{vwip}
(v^*vww^*)(v^*vww^*)
~=~
v^*vww^*.
}\eeq
Multiplying \eqref{vwip}  from the left with $v$ and from the right with $w$ gives back \eqref{vwpi}. So, both together show \ref{pi1} $\Leftrightarrow$ \ref{pi2}.

$v^*vww^*$ is a contraction. If it also is an idempotent, then \cite[Lemma B.1]{ShaSk23} asserts that it is a projection.\footnote{See Remark \ref{Hrem}.} So, \ref{pi2} $\Rightarrow$ \ref{pi3}.

If $v^*vww^*$ is a projection, then it is, in particular, self-adjoint. So,  \ref{pi3} $\Rightarrow$ \ref{pi4}.\qed

\bemp[Historical remark.~] \label{Hrem}
Partial isometries between Hilbert spaces have been studied since long.%
\footnote{
With this, we do not mean the occasions where they occur naturally such as basic things like \it{polar decomposition} or more advanced things like \it{Murray-von Neumann equivalence} in classification of operators algebras, but we mean papers that study partial isometries or subclasses of them for their own sake.
}
There is the amazing paper by Halmos and McLaughlin \cite{HaMLa63} who examine partial isometries \bf{on} Hilbert spaces under unitary equivalence; Halmos and Wallen \cite{HaWa70} determine the complete structure of discrete semigroups of partial isometries; Erdelyi (who has written quite a number of papers about partial isometries that, unfortunately, appear to be only little known) in \cite{Erd68} examined explicitly when partial isometries \bf{on} a Hilbert space have partial isometries as product.%
\footnote{
For isometries \bf{on} a Hilbert space the equivalence \eqref{pi1}$\Leftrightarrow$\eqref{pi2} in Lemma \ref{pilem} is, for instance, (including the simple proof) (i)$\Leftrightarrow$(iv) in \cite[Theorem 1]{Erd68}, and \eqref{pi1}$\Leftrightarrow$\eqref{pi4} in Lemma \ref{pilem} is, for instance, \cite[Lemma 2]{HaWa70}.
}
Apart from possibly removing some dimension restrictions (see Footnote \ref{onbetFN}), nothing much changes by passing to partial isometries \bf{between} Hilbert spaces. We guess that nothing about Lemma \ref{pilem} is new, except possibly the arrangement of the steps in the proof, and appealing -- to that goal -- to \cite[Lemma B.1]{ShaSk23}, which asserts, even for Hilbert modules, that a contractive right linear idempotent on a Hilbert module is a projection. The use of this lemma is \it{just for fun}; see the following remark.
\eemp

\bemp[About {\cite[Lemma B.1]{ShaSk23}}.~] \label{SSlemrem}
The proof of \cite[Lemma B.1]{ShaSk23} falls into two parts: First one shows that every contractive idempotent on a Hilbert \bf{space} is a projection. (We like to mention that this is true even for pre-Hilbert spaces.%
\footnote{
A contractive idempotent on a pre-Hilbert space $H$ extends continuously to a contractive idempotent on the completion $\ol{H}$. And after applying the lemma to this continuous extension, we obtain the statement for the original operator simply by (co)restriction to $H$, because the extended operator leaves $H$ invariant.
}%
) Then, the statement for Hilbert modules is boiled down to the Hilbert spaces case (using \bf{the tool} described in Footnote \ref{tpFN}) in a way very similar to how some proofs in Section \ref{modpiSEC} make use of the results for Hilbert spaces from the present section. Note that \cite[Lemma B.1]{ShaSk23} is, clearly, an analytic result. While the logic of the consecutiveness in the proof of Lemma \ref{pilem} would be a different one, we still like to mention that we can avoid using \cite[Lemma B.1]{ShaSk23} for proving Lemma \ref{pilem} taking into account that the idempotent $v^*vww^*$ to which we apply it has additional algebraic properties: 

We do have that \eqref{pi2} in Lemma \ref{pilem} ($v^*vww^*$ is an idempotent) implies \eqref{pi1} ($vw$ is a partial isometry). We will see in Footnote \ref{ivFN} that \eqref{pi4} ($v^*v$ and $ww^*$ commute) is equivalent to the property $ww^*v^*vw=v^*vw$, and, as mentioned in the proof of Lemma \ref{pilem}, \eqref{pi4} implies, obviously, \eqref{pi3} ($v^*vww^*$ is a projection). So, we are done showing that \eqref{pi2} implies \eqref{pi3}, if we show that \eqref{pi1} ($vww^*v^*vw=vw$) implies $ww^*v^*vw=v^*vw$ (equivalent to \eqref{pi4}). But this follow easily, showing (by multiple use of $vww^*v^*vw=vw$) that $(ww^*v^*vw-v^*vw)^*(ww^*v^*vw-v^*vw)=0$.
\eemp

\lf
The product $vw$ of two partial isometries is, in general, not again a partial isometry; but it is -- we used that already in the proof of Lemma \ref{pilem} -- always a contraction. We wish to turn the contraction into a partial isometry, following the idea of ``making the domain smaller without throwing away too much'' that worked already nicely for the composition of partially defined functions. Recall that for the partial isometry $v_f$ associated with an injective partially defined function $(f,D_f)$, the domain $D_f$ is encoded in the initial projection $v_f^*v_f$. We shall turn $vw$ into a partial isometry, making its ``domain'' smaller by multiplying $vw$ with a projection from the right -- the projection onto all those $x$ on which $vw$ acts isometrically.

\blem \label{clem}
Let $c$ be a contraction.
\begin{enumerate}
\item
Then $S=\bCB{x\colon\norm{cx}=\norm{x}}$ is a closed subspace.

\item
If $p$ is the projection onto $S$, then $(cp)^*(cp)=p$ and $p\le c^*c$.

\item
If $q$ is another projection such that $(cq)^*(cq)=q$, then $q\le p$.
\end{enumerate}
\elem

\proof
It is clear that $S$ is closed and invariant under multiplication with scalars in $\C$. Since $c^*c$ is a positive operator and a contraction, the same is true for the operator $a:=\U-c^*c$. Observe that
\beqn{
\norm{cx}=\norm{x}
~\Longleftrightarrow~
\AB{x,ax}=0
~\Longleftrightarrow~
\sqrt{a}x=0
~\Longleftrightarrow~
ax=0.
}\eeqn 
(From $ax=0$, we have $\AB{x,ax}=0$, hence, $\sqrt{a}x=0$.) So $S$ is a closed subspace -- the eigenspace of $a$ to the eigenvalue $0$, or, equivalently, the eigenspace of $c^*c$ to the eigenvalue $1$.

By the last property we have $c^*cp=p=pc^*c$, so $(cp)^*(cp)=pc^*cp=p$ and $p=pc^*cp\le pc^*cp+(\U-p)c^*c(\U-p)=c^*c$.

For $x=qx$, we find $\norm{cx}^2=\norm{cqx}^2=\AB{x,qc^*cqx}=\norm{qx}^2=\norm{x}^2$, so $x\in S$, hence $q\le p$.\qed

\lf
Note that by Lemma \ref{pilem}, $cp$ is, in  particular, a partial isometry with initial projection $p$. By definition of $p$, the partial isometry $cp$ coincides with $c$ on $S$. Henceforth, denote the projection $p$ constructed from a contraction $c$ according to Lemma \ref{pilem} by $p_c$.

Note that if $c=v$ is a partial isometry, then we have $p_v=v^*v$, the initial projection of $v$.

\bdefi \label{pisivdefi}
For a contraction $c$ we call $cp_c$ the \hl{partial isometry contained in $c$}.
\edefi

For the partial isometry $vwp_{v,w}$ contained in the contraction $vw$, we get:

\bthm \label{cthm}
For two partial isometries $v$ and $w$ there is a unique maximal projection $p_{v,w}$ such that $(vwp_{v,w})^*(vwp_{v,w})=p_{v,w}$.

Moreover, $p_{v,w}\le w^*w$, also $wp_{v,w}$ is a partial isometry, and $(wp_{v,w})(wp_{v,w})^*\le v^*v$.
\ethm

\proof
We apply Lemma \ref{clem} to the contraction $vw$ to get $p_{v,w}$. Moreover, $p_{v,w}\le(vw^*)(vw)=w^*v^*vw\le w^*w$. For every projection $p\le w^*w$ for a partial isometry $w$, we have $(wp)^*(wp)=p$ so that $wp$, especially our $wp_{v,w}$, is a partial isometry, too. Since, $wp_{v,w}$ is a partial isometry, $q:=(wp_{v,w})(wp_{v,w})^*$ is a projection. Since
\bmun{
~~~~~~
(vq)^*(vq)
~=~
\Bfam{(vwp_{v,w})(wp_{v,w})^*}^*\Bfam{(vwp_{v,w})(wp_{v,w})^*}
\\
~\le~
(wp_{v,w})(wp_{v,w})^*
~=~q,
~~~~~~
}\emun
applying Lemma \ref{clem} to the contraction $v$, we see $q\le v^*v$.\qed

\bob
For the partial isometries $v_f$ and $w_g$ coming from  injective partially defined functions, $p_{v_fv_g}$ is the initial projection of the isometry $v_fv_g=v_{f\circ g}$. It, therefore, reflects the making smaller of the domain of $g$ reflected by $w^*w$.
Also, the range projection $(wp_{v,w})(wp_{v,w})^*$ of the partial isometry $wp_{v,w}$ (representing, for functions, the restriction of $g$ to $D_{f\circ g}$), is dominated by the initial projection $v^*v$ (corresponding, for functions, to the fact that $g$ restricted to $D_{f\circ g}$ maps into the domain $D_f$ of $f$).
\eob

\bthm \label{cathm}
The composition $v\cdot w:=vwp_{v,w}$ among partial isometries that associates with two (composable) partial isometries the partial isometry contained in $vw$ is associative.
\ethm

\proof
In fact, we will show that both $(u\cdot v)\cdot w$ and $u\cdot(v\cdot w)$ are equal to $uvwp$, where $p$ is the projection in Lemma \ref{clem} constructed for the contraction $uvw$.

By the last part of Theorem \ref{cthm}, in the computation of $p_{uvp_{u,v},w}$ by applying Lemma \ref{clem} to the contraction $uvp_{u,v}w$ nothing changes if in the contraction $uvp_{u,v}w$ we leave out $p_{u,v}$, because $w$ maps $p_{uvp_{u,v},w}x$ into the range of $p_{u,v}$. Likewise,  in the computation of $p_{u,vwp_{v,w}}$ by applying Lemma \ref{clem} to the contraction $uvwp_{v,w}$ nothing changes if in the contraction $uvwp_{v,w}$ we leave out $p_{v,w}$, because $p_{u,vwp_{v,w}}$ must be not bigger than $p_{v,w}$. So in either case, we get $uvwp$ where $p$ is as stated in the beginning.\qed

\brem
From Theorem \ref{cathm} it follows that the class of Hilbert spaces with the partial isometries as morphisms forms a category under the composition ``$\cdot$'' of morphisms.
\erem

\newpage

\section{About composing partial isometries between Hilbert modules} \label{modpiSEC}

We now wish to push forward the considerations of the preceding section to partial isometries between Hilbert modules.%
\footnote{
\bf{For reader's convenience:} Recall that for a \nbd{C^*}algebra $\cB$, a Hilbert \nbd{\cB}module is a (right, of course) \nbd{\cB}module $E$ with a sesquilinear map $\AB{\bullet,\bullet}\colon E\times E\rightarrow\cB$ (the \hl{inner product)} satisfying:
\begin{itemize}
\item
$\AB{x,x}\ge 0$. (This plus sesquilinearity implies already symmetry $\AB{x,y}=\AB{y,x}^*$.)

\item
$\AB{x,yb}=\AB{x,y}b$. (The first two together imply already the (\nbd{C^*}algebra-valued!) Cauchy-Schwarz inequality
\beqn{
\AB{x,y}\AB{y,x}
~\le~
\AB{x,x}\norm{\AB{y,y}},
}\eeqn
from which it follows in the usual way that $\norm{x}:=\sqrt{\AB{x,x}}$ is a seminorm, that the kernel of this seminorm may be divided out, and that the normed space may be completed to a Banchch space inheriting the same structures.)

\item
$\AB{x,x}\Rightarrow x=0$.

\item
$E$ is complete in the norm $\norm{\bullet}$.
\end{itemize}
(See, for instance, Lance \cite{Lan95} or Skeide \cite{Ske01} for details.) By $\sB^r(E,F)$ we denote the bounded right linear maps between two Hilbert \nbd{\cB}modules; by $\sB^a(E,F)$ we denote the adjointable (see below) maps.
}
 Since we will no longer have available the luxury to be able to always write a closed submodule $E_0\subset E$ as $pE$ for a suitable projection, we will have to stick to the submodules themselves. And inequalities such as $c^*c\le d^*d$ will, in absence of an \it{adjoint} (see below), be occurring in the form $\AB{cx,cx}\le\AB{dx,dx}$.%
 \footnote{
 Recall that an operator $a\in\sB^r(E)$ is positive in the \nbd{C^*}algebra $\sB^a(E)\subset\sB^r(E)$ if and only if $\AB{x,ax}\ge0$ for all $x\in E$; see Paschke \cite[Proposition 6.1]{Pas73} or Lance \cite[Lemma 4.1]{Lan95}.
}

As we said, bounded right linear maps $a\colon E\rightarrow F$ between Hilbert \nbd{\cB}modules need not necessarily possess an \hl{adjoint}, that is, a map $a^*\colon F\rightarrow E$ such that $\AB{ax,y}=\AB{x,a^*y}$ for all $x\in E,y\in F$.%
\footnote{
Note that for any map $a\colon E\rightarrow F$ between Hilbert \nbd{\cB}modules, existence of an adjoint guarantees that $a$ is right linear and bounded. (If the modules fulfill everything required from a Hilbert module but possibly completeness, then an adjointable $a$ is right linear but not necessarily bounded.) Of course, $a^*$ is unique and $(a^*)^*=a$.
}%
\footnote{
All geometric type problems that Hilbert modules do in general pose, while Hilbert spaces do not, have a common root in the fact that contrary to Hilbert spaces, Hilbert modules need not be self-dual: Not every bounded right linear map $\Phi\colon E\rightarrow\cB$ need have the form $x\mapsto\AB{y,x}$ for a suitable $y\in E$. A related fact, requiring a really involved counter example to illustrate it, is that the extension of the \nbd{0}map on a submodule $F\subset E$ with zero-complement in $E$ to all of $E$ need not be unique; see Kaad and Skeide \cite{KaSk23}. (Note added: We feel the urge to mention the two papers, one by Manuilov  \cite{Man24} (first on arXiv) and another one by Frank \cite{FraM24} (second on arXiv), in which Manuilov \cite[Section 5]{Man24} constructs a counter example similar to \cite{KaSk23} with a Hilbert module over a monotone complete \nbd{C^*}algebra, while \cite[Theorem 4.4]{FraM24} asserts that precisely that would not be possible.) A subcategory of Hilbert modules that is free of all these problems is formed by von Neumann modules; see Skeide \cite{Ske00b,Ske01,Ske06b,Ske22b}.
}
This raises several questions such as adjointability of isometries (they need not be adjointable) and, worse, if things such as coisometries or partial isometries can be defined without (as the definition in the Hilbert space case in the preceding section suggests) assuming explicitly that they are adjointable.

A map $v\colon E\rightarrow F$ is usually called an \hl{isometry} if it preserves inner products. (These are exactly the embeddings of the structure.) It is easy to check that an isometry is right linear; of course, it is a contraction.%
\footnote{
One may show that a right linear map that is norm-preserving also preserves inner products. (See Lance \cite[Theorem 3.5]{Lan95}; Lance's lemma is actually about unitaries, but for the first part of its  proof it is not necessary to require that the map in question be surjective.)
}
But an isometry need not be adjointable.%
\footnote{
In Bhat and Skeide \cite[Section 2]{BhSk15} there are plenty of non-adjointable isometries -- even semigroups of such.
}
In fact, an isometry  is adjointable if and only if there exists a projections onto its range; see \cite[Proposition 1.5.13]{Ske01}. \hl{Projections} can be defined as maps $p\colon E\rightarrow E$ satisfying $\AB{px,py}=\AB{x,py}$.%
\footnote{
Examining $\abs{p(x+x'b)-px-(px')b}^2$ (where the algebra valued length is $\abs{y}:=\sqrt{\AB{y,y}}$), we see that $p$ is right linear. And by Cauchy-Schwarz inequality, $\norm{px}^2=\norm{\AB{x,px}}\le\norm{x}\norm{px}$, so $p$ is a contraction.
}
This is equivalent to requiring $p$ adjointable and $p^*p=p$. A \hl{unitary} is a surjective isometry. Of course, $u^{-1}=u^*$. But how do we define coisometries and partial isometries without requiring them \it{a priori} to be adjointable? (In fact, our coisometries will turn out to be adjointable, while our partial isometries need not be.)

We follow the definition proposed in Shalit and Skeide \cite[Remark 5.10]{ShaSk23}. The proposal exposed there, starts from the idea that the definition of coisometry should be a co-notation of the definition of isometry. For that goal it is necessary to give a different (but equivalent) definition of isometry. We quote from \cite[Remark 5.10]{ShaSk23}: ``An \hl{isometry} is [...] a contractive right linear map $E\rightarrow F$ between Hilbert modules for which there exists a submodule $F'$ of $F$ such that the map corestricts to a unitary $E\rightarrow F'$; a \hl{coisometry} is a contractive right linear map $E\rightarrow F$ between Hilbert modules for which there exists a submodule $E'$ of $E$ such that the map restricts to a unitary $E'\rightarrow F$.''

\bob
Clearly, the two definitions of isometry we stated are equivalent; in fact, the second is obviously stronger than the first, but checking just that inner products are preserved, is already enough to get right linearity and contractivity. Of course, the submodule $F'$ has no choice but to be the range of the isometry. Also the submodule $E'$ in the definition of coisometry is unique; but this can be seen easily only after showing in \cite[Corollary B.2 of Lemma B.1]{ShaSk23} that (unlike isometries) coisometries always have an adjoint. On the other hand, after knowing that a coisometry has an adjoint, we easily see that our definition of coisometry is equivalent to the definition for Hilbert spaces in Section \ref{picSEC}. (In fact, once a coisometry $v$ has an adjoint, it is clear that $v^*$ is an (adjointable!) isometry. And once this is clear, also the submodule $E'$ turns to be the range of $v^*$ (and the complement of the kernel of $v$).)
\eob

\lf
Finally, also for partial isometries we quote from \cite[Remark 5.10]{ShaSk23}:

\bemp[Definition \cite{ShaSk23}.]\hspace{-1ex}%
\footnote{
Halmos and McLaughlin \cite{HaMLa63} propose as definition to call a \it{partial isometry} a map that restricts to an isometry on the complement of its kernel. This is one of the thousands of equivalent ways to say what a partial isometry between Hilber spaces is, and it is, admittedly, one of the few that do not involve explicitly the adjoint of the map in question. However, for Hilbert modules it does not work. Indeed, by Kaad and Skeide \cite{KaSk23} there exist a Hilbert module, a submodule that has zero complement, and a non-zero bounded right linear map that vanishes on the submodule; according to the definition in \cite{HaMLa63}, such map would be a partial isometry. (In fact, such map need not even be a contraction as any of its multiples does satisfy that definition, too.)
}
~
A map $v\colon E\rightarrow F$ between Hilbert modules $E$ and $F$ is a \hl{partial isometry} if it corestricts to a coisometry $E\rightarrow vE$.
\eemp

Clearly, all operators we were talking about before (projections, isometries, coisometries, hence unitaries) are partial isometries.

Note that for a partial isometry $v$, the coisometry $v_s\colon E\rightarrow vE$ obtained by corestriction (the index ``s'' stands for ``surjective version'') does have an adjoint $v_s^*\colon vE\rightarrow E$ (which is an adjointable isometry). Therefore, a partial isometry $v$ does have an \hl{initial projection}, namely, $\pi_v:=v_s^*v_s$ (satisfying $v\pi_v=v$ and $(\U-\pi_v)E=\ker v$ and, obviously, being determined by these two properties).

Also, $v$ is adjointable if and only if $vE$ is complemented in $F$. (If $v$ is adjointable, then $vv^*$ is a projection onto $vE$; if $vE$ is complemented, then there is a projection $p\colon F\rightarrow F$ onto $vE$ and, denoting by $p_s\colon F\rightarrow vE$ the corestriction to $vE$ of $p$, it is a recommended exercise to check that $v_s^*p_s$ is an adjoint of $v$.)

We see, partial isometries inherit the good behavior of coisometries on the side of their domain (decomposing it into a direct sum of the range of their initial projection and their kernel), but they inherit the (potentially) bad behavior of the range of isometries (being not necessarily complemented).

\bob \label{picob}
Let us state some more facts about partial isometries. In the sequel, $w$ denotes a partial isometry $D\rightarrow E$ and $v$ denotes a partial isometry $E\rightarrow F$ (so that $vw\colon D\rightarrow  F$ makes sense).

\begin{enumerate}
\item \label{pic1}
By definition, the corestriction of $v$ as map from $E$ onto $vE$ is a coisometry. The restriction to a map from $\pi_vE$ to $F$ is an isometry. The (co)restriction to a map from $\pi_vE$ to $vE$ is a unitary.

\item \label{pic2}
If $v$ and $w$ are isometries, then $vw$ is an isometry. (Recall that for checking that it is enough to check the standard definition, that is, to check that $vw$ is inner product preserving. See also Remark \ref{injrem})

\item \label{pic3}
If $v$ and $w$ are coisometries, then $vw$ is a coisometry. (Recall that a map is a coisometry if and only it is adjointable and its adjoint is an isometry.)

\item \label{pic4}
If $v$ is an isometry, then $vw$ is a partial isometry. (The corestriction of $vw$ to $vwD$ restricts to a unitary on $\pi_wD$.)

\item \label{pic5}
If $w$ is a coisometry, then $vw$ is a partial isometry. (Since $w$ is a coisometry onto $E$, we may consider $\pi_vw$ a coisometry onto $\pi_vE$ and the restriction $v_v$ of $v$ to $\pi_vE$ is an isometry. Then the statement follows from $vw=v_vp_vw$ and the preceding points.)

\item
Especially, by \eqref{pic4} and \eqref{pic5}, if $v$ is an isometry and $w$ is a coisometry, then $vw$ is a partial isometry. But also (as well known for Hilbert spaces) the converse is true: Every partial isometry can be written as the product of an isometry and a coisometry. (Indeed, if $v$ is a partial isometry, then it can be written as the product of $v$ restricted to $\pi_vE$ -- an isometry -- and the corestriction of $\pi_v$ to $\pi_vE$ -- a coisometry.)%
\footnote{ \label{onbetFN}
Note that we cannot, in general, speak about partial isometries \bf{on} a Hilbert module (or space), but we have to speak about partial isometries between Hilbert modules (or spaces). For instance, already for a Hilbert space $H$, for any product of an isometry $v\in \sB(H)$ and a coisometry $w\in\sB(H)$ the dimensions of both $(vw)(vw)^*H=vv^*H$ and $(vw)^*(vw)H=w^*wH=w^*H$ cannot be smaller than $\dim H$.
}
\end{enumerate}
\eob

\lf\noindent
If both $v$ and $w$ are adjointable, then (especially, since \cite[Lemma B.1]{ShaSk23} also holds for Hilbert modules) nothing changes in the proof of Lemma \ref{pilem}. In the general case, $\pi_v=v_s^*v_s$ is a good substitute for $v^*v$, while we have no substitute for $ww^*$.  So there is not much we can do about the properties stated under \eqref{pi2} and \eqref{pi3} if $w$ is not adjointable. But if $v$ and $w$ should be adjointable, then \eqref{pi4} can be phrased equivalently as the condition that $\pi_v$ has to leave the range of $w$ invariant.%
\footnote{ \label{ivFN}
In the adjointable case, clearly, commutation of the projection $\pi_v=v^*v$ and projection $ww^*$ onto the range $wD$ of $D$ implies invariance of the latter under $\pi_v$. Conversely, invariance of $wD$ under $\pi_v$ is the same as $v^*vw=ww^*v^*vw$, so, $v^*vww^*=ww^*v^*vww^*$. Taking into account that the right-hand side is self-adjoint, we conclude $v^*vww^*=(v^*vww^*)^*=ww^*v^*v$.
}
And in this form, the condition makes sense also if $w$ is not adjointable.

\blem
Let $v\colon E\rightarrow F$ and $w\colon D\rightarrow E$ be partial isometries between Hilbert \nbd{\cB}modules. Then $vw$ is a partial isometry if and only if $\pi_v$ leaves $wD$ invariant, that is, for every $x\in D$ there is $y\in D$ such that $\pi_vwx=wy$.
\elem

\proof
The condition is sufficient: Indeed, $\pi_v$, sending $wD$ into $wD$ (co)restricts to a projection in $\pi_v^w\in\sB^a(wD)$, and $w$ corestricts to a(n adjointable!) coisometry $w_s$ onto $wD$. So, $\pi_v^ww_s$ is an adjointable partial isometry in $\sB^a(D,wD)$ with range contained in $\pi_vE\cap wD$.  Considering $\pi_v^ww_s$ as map into $\pi_vE$, we get a (not necessarily adjointable) map $w^p\in\sB^r(D,\pi_vE)$, obviously still a partial isometry (it corestricts to a coisometry onto $\pi_vE\cap wD\subset \pi_vE$). Now, $v$, on $\pi_vE$, acts isometrically, and since the composition of an isometry with a partial isometry is a partial isometry (Observation \ref{picob}\eqref{pic4}), the composition $(v\upharpoonright \pi_vE)w^p=vw$ is a partial isometry.

The condition is also necessary: Indeed, suppose $vw$ is a partial isometry. Since $v_sE=vE\supset vwD$, also $v_sw$ is a partial isometry. Composing it with the isometry $v_s^*$, we see that $\pi_vw=v_s^*v_sw$ is a partial isometry. Also, since $w\pi_w=w$ and $\pi_{vw}\le\pi_w$ (this follows from $\AB{vwx,vwx}\le\AB{wx,wx}$ which by Paschke \cite[Theorem 2.8]{Pas73} holds for every contraction $v$), we may restrict to $\pi_w D$ assuming, thus, that $w$ is an isometry. So we have reduced the problem to the case of an isometry $w\colon D\rightarrow E$ and a projection $p\in\sB^a(E)$ such that $pw$ is a partial isometry. In particular, there is the projection $\pi_{pw}$ such that $pw$ is an isometry on $\pi_{pw}D$ and $pw$ is $0$ on $(\U-\pi_{pw})D$. For $x\in\pi_{pw}D$ we find
\beqn{
\AB{wx-pwx,wx-pwx}
~=~
\AB{wx,wx}-\AB{pwx,pwx}
~=~
0,
}\eeqn
so $pwx=wx\in wD$. For $y\in(\U-\pi_{pw})D$ we find $pwy=0$. A general element of $D$ has the form $x+y$, so a general element of $wD$ has the form $w(x+y)$. Since $pw(x+y)=wx$, we see that $p$ leaves $wD$ invariant.\qed

\lf
We now pass to the module analogue of Lemma \ref{clem}. As stated in the introduction to this section, we do, in general, not have projections onto closed submodules, so we have to stick with the submodules. But the construction of the desired submodule is similar to Lemma \ref{clem} and the proof here makes use of it.

\bthm \label{cMthm}
Let $c\in\sB^r(E,F)$ be a contraction. Then
\beqn{
P_c
~:=~
\bCB{x\colon\AB{cx,cx}=\AB{x,x}}
}\eeqn
is the unique maximal closed submodule on which $c$ acts isometrically.
\ethm

\proof
As in Lemma \ref{clem}, $P_c$ is clearly closed and invariant under right multiplication. If it is also invariant under addition, then, by polarization, $c$ acts isometrically on $P_c$ and there are no other elements in $E$ on which $c$ could do that.

In a way very similar to the proof of \cite[Lemma B.1]{ShaSk23}, by the following tool -- we say \bf{the tool} -- we reduce the proof that $P_c$ is closed under addition in the module case to the case of Hilbert spaces in Lemma \ref{clem}. Assume that $\cB$ is acting as a concrete operator algebra $\cB\subset\sB(G)$ (nondegenerately) on a Hilbert space $G$. As explained in Skeide \cite{Ske00b} (see also \cite{Ske22b}), we may construct a Hilbert space $H$ and for every $x\in E$ an operator $L_x\in\sB(G,H)$ such that $\AB{x,y}=L_x^*L_y$ ($x,y\in E$). Necessarily, $L_{xb}=L_xb$ and every operator $c\in\sB^r(E)$ gives rise to an operator $C\in\sB(H)$ satisfying $CL_x=L_{cx}$ and $\norm{C}=\norm{c}$. It follows that $\AB{L_xg,L_yg'}=\AB{g,\AB{L_x,L_y}g'}$.%
\footnote{ \label{tpFN}
For those who know the (internal) tensor product: $H=E\odot G$, the tensor product of the Hilbert \nbd{\cB}module with the correspondence $G$ from $\cB$ to $\C$ over $\cB$ with $L_x=x\odot\id_G\colon g\mapsto x\odot g$ and $C=c\odot\id_G$. While $\norm{C}=\norm{c}$ for adjointable $c$ follows easily from automatic contractivity (isometricity) of (faithful) homomorphisms between \nbd{C^*}algebras, for non-adjointable $c$ this is harder to prove; see \cite[Appendix]{Ske22b}.
}
The same goes  for our contraction $c\in\sB^r(E,F)$, when $F\subset\sB(G,K)$ for another Hilbert space $K$. Applying Lemma \ref{clem} to $H$ and the contraction $C$, we find $L_xg\in S\subset H$ for all $x\in P_c$ and $g\in G$. Since $S$ is a subspace, it follows that $L_{x+y}(g+g')\in S$ for all $x,y\in P_c$ and $g,g'\in G$. By polarization, we conclude
\bmun{
\BAB{g,\bAB{c(x+y),c(x+y)}g'}
~=~
\AB{CL_{x+y}g,CL_{x+y}g'}
\\
~=~
\AB{L_{x+y}g,L_{x+y}g'}
~=~
\BAB{g,\bAB{x+y,x+y}g'}
}\emun
for all $x,y\in P_c$ and $g,g'\in G$, so $\bAB{c(x+y),c(x+y)}=\bAB{x+y,x+y}$, so $x+y\in P_c$.\qed

\lf\lf
So for a contraction (such as $vw$) we do find a maximal submodule on which it behaves isometric. But this is as far as we go in our attempts to possibly find a partial isometry that coincides with $c$ on $P_c$ and tells us enough about $P_c$.

\bex
Let $E=\cB=C\SB{0,2}$ so that $\sB^r(E)=\sB^a(E)=\cB$ (having act $a\in\sB^r(E)=\cB$ on $E=\cB$ by multiplication (from the left, obviously)). Now since $\cB$ is commutative, the only partial isometries are projections; since the only projection in $\cB$ is $\U$, there are no partial isometries but those of the form $\U e^{\imath\vp}$ ($\vp\in\R$).

Let $c\in C_0\LO{0,2}\subset\cB$ such that $0<c(t)<1$ for all $t\in(0,1)$ and $c(t)=1$ for all $t\in\SB{1,2}$. Then $c$ acts as a contraction (even adjointable) on $E$ and, clearly,
\beqn{
P_c
~=~
\bCB{b\in\cB\colon b(t)=0\text{~for all~}t\in\SB{0,1}}.
}\eeqn
But since $\cB^r(E)$ contains up to a phase factor only one partial isometry and the only one among all those that coincides on $P_c$ with $c$ is $\U$, that partial isometry about $P_c$ tells us absolutely nothing.%
\footnote{
This example is, clearly, promotable to other Hilbert modules over that commutative $\cB$ or more general \nbd{C^*}algebras with a center that has a nontrivial connected part. We propose for future work to examine the question for modules (preferably, full) over $\cB$ (preferably, unital) with trivial center (distinguishing also between simple and not).
}
\eex

\lf
While it occurs difficult to realize $c$ in that example as a product $c=vw$ of partial isometries $v$ and $w$, we show now that the hope to be able to associate to a contraction $c$ ($=vw$ or general) a partial isometry characterized by properties similar to those in Lemma \ref{clem} (but phrased without assuming existence of adjoints or projections) is vane, unless $P_c$ is complemented or, equivalently, unless there exists a projection $p_c$ onto $P_c$.

Recall that Lemma \ref{clem} grants for any contraction $c$ between Hilbert spaces existence of a partial isometry $v=cp_c$ (where $p_c$ is the projection onto $S=P_c$; see the discussion leading to Definition \ref{pisivdefi}) that coincides with $c$ for elements in $P_c$ and is such that $p_c=\pi_v=v_s^*v_s\le c^*c$. Since, for $c\colon E\rightarrow F$ between modules, we do not necessarily have an adjoint of $c$, we write this inequality as
\beq{ \label{wdomeq}
\AB{x,\pi_vx}
~\le~
\AB{cx,cx}
}\eeq
for all $x\in E$ (which, should $c^*$ exist, is equivalent to $\pi_v\le c^*c$).

\bdefi
Let $c\colon E\rightarrow F$ be a right linear contraction. If $v\colon E\rightarrow F$ is a partial isometry satisfying \eqref{wdomeq} and $vx=cx$ for all $x\in P_c$, then we call $v$ \hl{the} ($v$ turns out to be unique) \hl{partial isometry contained in $c$}.
\edefi

\bthm \label{univthm}
Let $c\colon E\rightarrow F$ be a right linear contraction.
\begin{enumerate}
\item \label{uv1}
If there exists a partial isometry $v\colon E\rightarrow F$ contained in $c$, then $v$ is unique and $\pi_v$ is a projection onto $P_c$.

\item\label{uv2}
If $p_c$ is a projection onto $P_c$, then $v=cp_c$ is a partial isometry, the partial isometry contained in $c$.
\end{enumerate}
Summing up, the partial isometry contained in $c$ exists if and only if $P_c$ is complemented.
\ethm

\proof
Part \eqref{uv1}: Since $c$ is a right linear contraction, we have $\AB{x,x}\ge\AB{cx,cx}$ for all $x\in E$. For $x\in\pi_vE$, we also have $\AB{x,x}=\AB{x,\pi_vx}\le\AB{cx,cx}$ by \eqref{wdomeq}, so $\AB{x,x}=\AB{cx,cx}$, hence, $\pi_vE\subset P_c$. On the other hand, since $vx=cx$ for all $x\in P_c$, we see $\AB{x,x}=\AB{cx,cx}=\AB{vx,vx}=\AB{x,\pi_vx}$, so $x\in\pi_vE$, hence $P_c\subset\pi_vE$. Since $\pi_v$ is unique, and $v$ is determined by what it does on $\pi_vE$, also $v$ is unique.

Part \eqref{uv2}: Since $c$ acts isometrically on $P_c=p_cE$, clearly, $v$ is a coisometry onto $vE$, that is, $v$ is a partial isometry and $\pi_v=p_c$. By definition, it satisfies $vx=cx$ for $x\in P_c=p_cE$. Already in the proof of Lemma \ref{clem}, where $c$ was adjointable, for the line that proves $p_c\le c^*c$ (which, in the adjointable case, is equivalent to \eqref{wdomeq}), it was crucial to know that $p_c$ is the projection onto an eigen space of $c^*c$ and, therefore, commutes with $c^*c$ so that, in particular, $(\U-p_c)c^*cp_c=0$. This condition is the same as
\beq{ \label{wcommeq}
\AB{c(\U-p_c)y,cp_cx}
~=~
0
}\eeq
for all $x,y\in E$. And in this form, it makes sense also for non-adjointable $c$ and still proves \eqref{wdomeq} also in that case. To see that \eqref{wcommeq} is true, we proceed as in the proof of Theorem \ref{cMthm}, and lift everything via \bf{the tool} to operators between the Hilbert spaces $H$ ($:=E\odot G$) and $K$ ($:=F\odot G$); see Footnote \ref{tpFN} for the terms in brackets. We get a contraction $C$ ($:=c\odot\id_G$) and a partial isometry $V$ ($:=V\odot\id_G$) from $H$ to $K$, and the initial projection $V^*V$ of $V$ coincides with $\hat{P}_c$ ($:=p_c\odot\id_G$). We do not care if $V$ is the maximal partial isometry $Cp$ constructed in Lemma \ref{clem} from the maximal projection $p$ onto $S\subset H$. (It probably is, but we need not know.) What we care about is that $\hat{P}_c$ is, clearly, a subprojection of that $p$ ($C$ acts isometrically on $\hat{P}_cH$, so $\hat{P}_cH\subset S=pH$), and if $p$ commutes with $C^*C$, then so does $P_c$. And from $(\U-\hat{P}_c)C^*C\hat{P}_c=0$ it follows, by taking matrix elements $\AB{g,\bullet g'}$, that also \eqref{wcommeq} holds.\qed

\lf
The following example illustrates that $P_{vw}$ need not be complemented.

\bex
Let $E$ be a Hilbert \nbd{\cB}module with a closed submodule $F\ne E$ such that $F^\perp=\zero$; in particular, $F$ is not complemented in $E$. On $E^2=\rtMatrix{E\\E}$, define the projection
\beqn{
p
~:=~
\frac{1}{2}\SMatrix{1&1\\1&1}
\colon
\SMatrix{x\\y}
~\longmapsto~
\frac{1}{2}\SMatrix{x+y\\x+y}.
}\eeqn
Then $p\rtMatrix{E\\F}=\bCB{\rtMatrix{x\\x}\colon x\in E}$. For the contraction $c$ obtained by composing the canonical injection $\rtMatrix{E\\F}\rightarrow\rtMatrix{E\\E}$ (a partial isometry) with the projection $p$ (also a partial isometry), we see that $P_c=\bCB{\rtMatrix{y\\y}\colon y\in F}$. ($p$ acts isometrically only on elements of the form $\rtMatrix{x\\x}$ and in the range of the canonical injection the lower component has to be in $F$.) Clearly, $P_c$ is not complemented in $\rtMatrix{E\\F}$.
\eex

\lf
In conclusion, for Hilbert modules we cannot propose the composition $v\cdot w$ assigning as for Hilbert spaces to two (composable) partial isometries a partial isometry on the whole space. We can, however, return to what we did in Section \ref{introSEC} for functions, defining them only partially, introducing the class of \hl{partially defined isometries} $(v,D_v)$, meaning exactly the same as a partially defined function, just that the class of sets is restricted to Hilbert \nbd{\cB}modules, the partial domain $D_v$ has to be a closed submodule (instead of just any subset), and $v$ has to act on $D_v$ as an isometry. Then, clearly, given a partially defined isometry $(w,D_w)$ from $D$ to $E$ and a partially defined isometry $(v,D_v)$ from $E$ to $F$, then $(v\circ w,D_{v\circ w})$ as defined in Equation \eqref{pcdef} in the introduction is a partially defined isometry from $D$ to $F$. (In fact, if $x\in D_{v\circ w}$, then $w$ sends $x$ (isometrically!) into the domain $D_v$ of $v$, from which $v$ sends $wx$ further (isometrically!) into $F$.)

\brem
Clearly, the class of Hilbert \nbd{\cB}modules with the class of partially defined isometries as morphisms is a category, in fact, a subcategory of the class of sets with the (injective) partially defined functions.
\erem

Let us summarize: By Theorem \ref{cMthm}, to any right linear contraction $c$ we construct the unique maximal submodule $P_c$ of its domain such that the restriction of $c$ to $P_c$ is an isometry. We denote the corresponding partially defined isometry by $(c,P_c)$.

\bdefi
Let $c\colon E\rightarrow F$ be a right linear contraction. Then we call $(c,P_c)$ \hl{the partially defined isometry contained in $c$}.
\edefi

Note that, unlike the \it{partial isometry} contained in $c$, the \it{partially defined isometry} contained in $c$ always exists.

Let us close with the relationship between the (not necessarily partially isometric) usual product $vw$ of two partial isometries $v$ and $w$, and their composition (always a partially defined isometry!) when viewing them as partially defined isometries $(v,\pi_vE)$ and $(w,\pi_wD)$. Without the obvious proof, we state:

\bprop
\hfill
$(vw,P_{vw})~=~(v,\pi_vE)\circ(w,\pi_wD)$.
\hfill\hfill\hfill
{~}
\eprop

\newpage

\vspace{1ex}\noindent
\bf{Acknowledgments.~}
I wish to thank Malte Gerhold and Orr Shalit for useful discussions, and Raja Bhat and Arup Pal for useful bibliographical hints.

\vspace{-2ex}
\newcommand{\Swap}[2]{#2#1}\newcommand{\Sort}[1]{}
\providecommand{\bysame}{\leavevmode\hbox to3em{\hrulefill}\thinspace}
\providecommand{\MR}{\relax\ifhmode\unskip\space\fi MR }
\providecommand{\MRhref}[2]{%
  \href{http://www.ams.org/mathscinet-getitem?mr=#1}{#2}
}
\providecommand{\href}[2]{#2}


\end{document}